\documentclass[a4paper,12p]{article}
\title{Constancy of curvature and conformal-projective flatness of statistical manifolds}
\author{\textrm{Chol Rim Min, Won Hak Ri, Kum Hyok Kwak } \\  
                  {\small\textit{Faculty of Mathematics, }}
                        {\small\textit{Kim Il Sung University,  Democratic Peoples' Republic of Korea}}\\   
}
    
\date{June, 2016}  
\usepackage{latexsym}   
\usepackage[pdftex]{graphicx}
\usepackage{amsmath}
\newtheorem{Thot}{Theorem}[section]
\newtheorem{Prop}{Proposition}[section] 
\newtheorem{Coro}{Corollary}[section]

\newtheorem{lemm}{Lemma}[section]
    
\begin{document}
\maketitle      
\begin{abstract}   
An identity of conformal-projective curvature tensor of a statistical manifold $(M, g, \nabla)$  is studied in this paper. The relation between the constancy of curvature and conformal-projective flatness of statistical manifolds is also discussed.
\\

{\small\textit{Keywords}: statistical manifold of constant curvature, conformal-projective flatness of statistical manifolds, conformal-projective curvature tensor}\\ 
\end{abstract}
\section{Introduction}
\quad 
Conformal-projective equivalence of statistical manifolds can be considered as a natural generalization of conformal equivalence of Riemannian metrics, which was introduced by Matsuzoe [4].

Conformal equivalence of Riemannian metrics and projective equivalence of affine connections are combined or generalized to lead up to conformal-projective equivalence or $\alpha-$conformal equivalence for statistical manifolds(See [2, 4, 5]).
Conformal-projective curvature tensor in a statistical manifold plays an important role as Weyl conformal curvature tensor does in Riemannian geometry.

Let $M$ be an $n-$dimensional manifold, $\nabla$ a torsion-free affine connection on $M$ , and $g$ a Riemannian metric on $M$ .

We denote by $\Gamma(E)$ the set of smooth sections of a vector bundle $E\to M$. So $\Gamma(TM)$ means the set of smooth vector field on $M$ and 
$\Gamma(TM^{(r,s)})$ means the set of tensor fields of type $(r,s)$ on $M$. The curvature tensor of  $\nabla$ is denoted by $R\in\Gamma(TM^{(1,3)})$.

Two statistical manifold $(M, g, \nabla)$ and $(M, \bar{g}, \bar{\nabla})$ are said to be conformally-projectively equivalent (or generalized conformal equivalent) if there exist two functions $\varphi, \psi\in C^{\infty}(M)$ satisfying that 
\begin{align}
 &\bar{g}(X,Y)=e^{\varphi+ \psi}g(X,Y)\nonumber\\
 &\bar{\nabla }_{X}Y=\nabla _{X}Y+d\varphi(X)Y+d\varphi(Y)X-g(X,Y)grad_{g}\psi\nonumber
 \end{align}
 for all $X,Y\in\Gamma(TM)$ (See [4]).
 
 Conformal-projective curvature tensor $W$ is defined by
 \begin{align}
  W(X,Y)Z &=R(X,Y)Z+\frac{1}{n(n-2)}\{Y[(n-1)Ric(X,Z)+\stackrel{*}{Ric}(X,Z)]-\nonumber\\
 & -X[(n-1)Ric(Y,Z)+\stackrel{*}{Ric}(Y,Z)]+[(n-1)\stackrel{*}{Ric}^\sharp (Y)+Ric^\sharp (Y)]g(X,Z)- \tag{1.1} \\
 & -[(n-1)\stackrel{*}{Ric}^\sharp (X)+Ric^\sharp (X)]g(Y,Z)\}+\frac{\sigma}{(n-1)(n-2)}[Xg(Y,Z)-Yg(X,Z)]\nonumber
 \end{align}
 and an  $n(\ge 4)-$dimensional statistical manifold $(M, g, \nabla)$ is conformally-projectively flat if and only if the conformal-projective curvature tensor vanishes everywhere on $M$ ,where
$R, Ric$ and $\sigma$ are curvature tensor field, Ricci tensor field and scalar curvature of $(M, g, \nabla)$ , respectively, and the Ricci operator $Ric^{\sharp}$ of $(M, g, \nabla)$ is the (1, 1)-tensor field determined by
\begin{equation*}
	g(Ric^{\sharp}(X),Y)=Ric(X,Y)
\end{equation*}
and the corresponding quantities for $\stackrel{*}\nabla$ which is a dual connection of $\nabla$ are denoted with $\ast$  (See [3]).

On the other hand, for any  $\alpha\in\mathbf{R}$ , two statistical manifolds $(M, g, \nabla)$ and $(M, \bar{g}, \bar{\nabla})$ are said to be $\alpha$ -conformally equivalent if there exists a function $\varphi\in C^{\infty}(M)$ satisfying that
\begin{align}
 &\bar{g}(X,Y)=e^{\varphi}g(X,Y)\nonumber\\
 &g(\bar{\nabla }_{X}Y,Z)=g(\nabla _{X}Y,Z)-\frac{1+\alpha}{2}d\varphi(Z)g(X,Y)+\frac{1-\alpha}{2}[d\varphi(X)g(Y,Z)+d\varphi(Y)g(X,Z)]\nonumber
 \end{align}
 for all $X,Y,Z\in\Gamma(TM)$ (See [2]).

It is easily verified that if two statistical manifolds  $(M, g, \nabla)$ and $(M, \bar{g}, \bar{\nabla})$ are 1-conformal equivalent then they are conformal-projective equivalent and so if a statistical manifold $(M, g, \nabla)$ is 1-conformal flat then it is conformal-projective flat.

We study some properties of conformal-projective curvature tensor of a statistical manifold and a sufficient condition for a statistical manifold to be conformal-projective flat in this paper.

In section 2, we show some properties of conformal-projective curvature tensor of a statistical manifold.

In section 3, we show that a statistical manifold of constant curvature is conformal-projective flat.
\section{Conformal-projective curvature tensor of a statistical manifold}
\quad 
In this section, we give an expression of conformal-projective curvature tensor $W$ of a statistical manifold $(M, g, \nabla)$ and show the relationship between conformal-projective curvature tensor $W$ of a statistical manifold $(M, g, \nabla)$ and one $\stackrel{*}{W}$ of a dual statistical manifold $(M, g, \stackrel{*}\nabla)$ .

We first give the following fact, which has been obtained independently by Zhang [6]. We quote the fact from his paper with a suitable modification for later use.
\begin{Prop}\label{pro2_1} $([6])$
 Let $\sigma$ and  $\stackrel{*}\sigma$ be the scalar curvature of a statistical manifold $(M, g, \nabla)$ and a dual statistical manifold $(M, g, \stackrel{*}\nabla)$ , respectively. Then we have \begin{equation}
 \sigma=\stackrel{*}\sigma \tag{2.1}
\end{equation}
\end{Prop}
\begin{lemm}\label{lemm2_1}
The conformal-projective curvature tensor $W$ of a statistical manifold $(M, g, \nabla)$ can be expressed as follows:
\begin{equation}
 W(X,Y)Z=R(X,Y)Z+YL(X,Z)-XL(Y,Z)+ \stackrel{*}{L}^{\sharp}(Y)g(X,Z)-\stackrel{*}{L}^{\sharp}(X)g(Y,Z)\tag{2.2}
\end{equation}
for all $X,Y,Z\in\Gamma(TM)$, where $L,\stackrel{*}{L}$ and $\stackrel{*}{L}^{\sharp}$ are tensor fields of type $(0, 2)$ and $(1, 1)$, respectively, given by
\end{lemm} 
\begin{align}
 &L(X,Y)=\frac{1}{n-2}\{\frac{1}{n}[(n-1)Ric(X,Y)+\stackrel{*}{Ric}(X,Y)]-\frac{\sigma}{2(n-1)}g(X,Y)\}\nonumber\\
 &\stackrel{*}{L}(X,Y)=\frac{1}{n-2}\{\frac{1}{n}[(n-1)\stackrel{*}{Ric}(X,Y)+Ric(X,Y)]-\frac{\stackrel{*}\sigma}{2(n-1)}g(X,Y)\}\nonumber\\
 &g(\stackrel{*}{L}^{\sharp}(X),Y)=\stackrel{*}{L}(X,Y)\nonumber
\end{align}
\textbf{Proof}. Using Eq. (2.1), we can express Eq. (1.1) as follows:
\begin{align}
W(X,Y)Z &=R(X,Y)Z+\frac{1}{n(n-2)}\{Y[(n-1)Ric(X,Z)+\stackrel{*}{Ric}(X,Z)]-\nonumber\\
& -X[(n-1)Ric(Y,Z)+\stackrel{*}{Ric}(Y,Z)]+[(n-1)\stackrel{*}{Ric}^{\sharp}(Y)+Ric^{\sharp}(Y)]g(X,Z)-\nonumber\\
& -[(n-1)\stackrel{*}{Ric}^{\sharp}(X)+Ric^{\sharp}(X)]g(Y,Z)\}+\frac{\sigma+\stackrel{*}\sigma}{2(n-1)(n-2)}[Xg(Y,Z)-Yg(X,Z)]\nonumber\\
& =R(X,Y)Z+\frac{1}{n-2}Y\{\frac{1}{n}Ric(Y,Z)+\stackrel{*}{Ric}(X,Z)]-\frac{\sigma}{2(n-1)}g(X,Z)-\nonumber\\
& -\frac{1}{n-2}X\{\frac{1}{n}[(n-1)Ric(Y,Z)+\stackrel{*}{Ric}(Y,Z)]-\frac{\sigma}{2(n-1)}g(Y,Z)+\nonumber\\
& +\frac{1}{n-2}\{\frac{1}{n}[(n-1)\stackrel{*}{Ric}^{\sharp}(Y)+Ric^{\sharp}(Y)]-\frac{\stackrel{*}\sigma}{2(n-1)}g(X,Z)+\nonumber\\
& +\frac{1}{n-2}\{\frac{1}{n}[(n-1)\stackrel{*}{Ric}^{\sharp}(X)+Ric^{\sharp}(X)]-\frac{\stackrel{*}\sigma}{2(n-1)}g(Y,Z)\nonumber
\end{align}
On the other hand, since $g(\stackrel{*}{Ric}^{\sharp}(X),Y)=\stackrel{*}{Ric}(X,Y)$ and $g(Ric^{\sharp}(X),Y)=Ric(X,Y)$ hold for all $X,Y\in\Gamma(TM)$ , we have
\begin{equation*}
g((n-1)\stackrel{*}{Ric}^{\sharp}(X)+Ric^{\sharp}(X)-\frac{\stackrel{*}\sigma}{2(n-1)}X,Y)=\stackrel{*}{L}(X,Y)
\end{equation*}
for all $X,Y\in\Gamma(TM)$.
So 
\begin{equation*}
\stackrel{*}{L}^{\sharp}(X)=(n-1)\stackrel{*}{Ric}^{\sharp}(X)+Ric^{\sharp}(X)-\frac{\stackrel{*}\sigma}{2(n-1)}X
\end{equation*}
holds for all $X\in\Gamma(TM)$.
Therefore we have
\begin{equation*}
W(X,Y)Z=R(X,Y)Z+YL(X,Z)-XL(Y,Z)+ \stackrel{*}{L}^{\sharp}(Y)g(X,Z)-\stackrel{*}{L}^{\sharp}(X)g(Y,Z)
\end{equation*}
for all $X,Y,Z\in\Gamma(TM)$. $\quad\Box$ \\
\begin{Thot}\label{theo2_1}
Let $W$ and $\stackrel{*}W$ be the conformal-projective curvature tensors of a statistical manifold $(M, g, \nabla)$ and a dual statistical manifold $(M, g, \stackrel{*}\nabla)$ , respectively. Then we have 
\begin{equation}
g(W(X,Y)Z,U)+g(\stackrel{*}{W}(X,Y)U,Z)=0\tag{2.3}
\end{equation}
for all $X,Y,Z,U\in\Gamma(TM)$.
\end{Thot}
\textbf{Proof}. From Eq. (2.2), we have
\begin{align}
g(W(X,Y)Z,U)&=g(R(X,Y)Z,U)+g(Y,U)L(X,Z)-g(X,U)L(Y,Z)+\nonumber\\
&+\stackrel{*}{L}(Y,U)g(X,Z)-\stackrel{*}{L}(X,U)g(Y,Z)\nonumber\\
g(\stackrel{*}{W}(X,Y)U,Z)&=g(\stackrel{*}{R}(X,Y)U,Z)+g(Y,Z)\stackrel{*}{L}(X,U)-g(X,Z)\stackrel{*}{L}(Y,U)+\nonumber\\
&+L(Y,Z)g(X,U)-L(X,Z)g(Y,U)\nonumber
\end{align}
for all $X,Y,Z,U\in\Gamma(TM)$. Since $g(R(X,Y)Z,U)+g(\stackrel{*}{R}(X,Y)U,Z)=0$ holds for all  $X,Y,Z,U\in\Gamma(TM)$, Eq. (2.3) holds.$\quad\Box$ \\

Eq. (2.3) shows that a statistical manifold $(M, g, \nabla)$ is conformally-projectively flat if an only if the dual statistical manifold $(M, g, \stackrel{*}\nabla)$ is conformally-projectively flat.
\section{Constancy of curvature and conformal-projective flatness of a statistical manifold}
\quad 
It is known that if an $n(\ge4)-$dimensional Riemannian manifold is of constant curvature, it is conformal flat and that an $n(\ge3)-$dimensional Riemannian manifold is projective flat if and only if it is of constant curvature.

The conformal-projective equivalence of a statistical manifold is a generalization of conformal equivalence and projective equivalence of a Riemannian manifold. So constancy of curvature of a statistical manifold has a close relationship to the conformal-projective equivalence of a statistical manifold.
It is shown that an  $n(\ge2)-$dimensional statistical manifold $(M, g, \nabla)$ is of constant curvature if and only if the tangent bundle $TM$ over $M$ with complete lift statistical structure $(g^{c}, \nabla^{c})$ is conformally-projectively flat by Hasegawa [1]. 

The following theorem shows that the relationship between constancy of curvature and conformal-projective flatness of a statistical manifold $(M, g, \nabla)$.
\begin{Thot}\label{theo3_1}
If an $n(\ge4)-$dimensional statistical manifold  $(M, g, \nabla)$ is of constant curvature, it is conformally-projectively flat.
\end{Thot}
\textbf{Proof}. From Eq. (2.2), we have
\begin{equation*}
W(X,Y)Z=R(X,Y)Z+YL(X,Z)-XL(Y,Z)+ \stackrel{*}{L}^{\sharp}(Y)g(X,Z)-\stackrel{*}{L}^{\sharp}(X)g(Y,Z)
\end{equation*}
for all $X,Y,Z\in\Gamma(TM)$.

Since a statistical manifold  $(M, g, \nabla)$ is of constant curvature,
\begin{equation*}
Ric=\stackrel{*}{Ric}
\end{equation*}
holds and so we have
\begin{equation*}
L(X,Y)=\stackrel{*}{L}(X,Y)=\frac{1}{n-2}\{Ric(X,Y)-\frac{\sigma}{2(n-1)}g(X,Y)\}
\end{equation*}
for all $X,Y\in\Gamma(TM)$. On the other hands, since
\begin{equation*}
R(X,Y)Z=K\{g(Y,Z)-g(X,Z)\}
\end{equation*}
holds for all $X,Y,Z\in\Gamma(TM)$, we have
\begin{equation*}
Ric(Y,Z)=tr\{X\mapsto R(X,Y)Z\}=(n-1)Kg(Y,Z)
\end{equation*}
for all $Y,Z\in\Gamma(TM)$. Since from the above equation,
\begin{equation*}
\sigma=tr_{g}\{(Y,Z)\mapsto Ric(Y,Z)\}=n(n-1)K
\end{equation*}
holds, we have
\begin{equation*}
L(X,Y)=\stackrel{*}{L}(X,Y)=\frac{K}{2}g(X,Y)
\end{equation*}
for all $X,Y\in\Gamma(TM)$ and so 
\begin{equation*}
L^\sharp(X)=\stackrel{*}{L}^{\sharp}(X)=\frac{K}{2}X
\end{equation*}
holds for all $X\in\Gamma(TM)$.
Therefore we have 
\begin{align}
W(X,Y)Z & =R(X,Y)Z+YL(X,Z)-XL(Y,Z)+ \stackrel{*}{L}^{\sharp}(Y)g(X,Z)-\stackrel{*}{L}^{\sharp}(X)g(Y,Z)\nonumber\\
& =K\{Xg(Y,Z)-Yg(X,Z)\}+\frac{K}{2}Yg(X,Z)-\frac{K}{2}Xg(Y,Z)+\nonumber\\
& +\frac{K}{2}Yg(X,Z)-\frac{K}{2}Xg(Y,Z)\nonumber\\
& =0\nonumber
\end{align}
for all $X,Y,Z\in\Gamma(TM)$ . So the proof is finished.$\quad\Box$ \\

If a statistical manifold is a self-dual statistical manifold , that is, a Riemannian manifold, conformal-projective flatness of a statistical manifold becomes to conformal flatness of a Riemannian manifold. So theorem 3.1 shows that if a Riemannian manifold is of constant curvature, it is conformal flat , which is well known in Riemannian geometry.

Consequently, theorem 3.1 generalizes the fact that an $n(\ge4)-$dimensional Riemannian manifold of constant curvature is conformal flat to case of a statistical manifold.

Theorem 3.1 and theorem in [1] give the following:
\begin{Coro}\label{coro2_1}
  Let  $(M,g,\nabla)$ be an $n(\ge4)-$dimensional statistical manifold. If the tangent bundle $(TM,g^c,\nabla^c)$ over $M$ with complete lift statistical structure is conformally-projectively flat, $(M,g,\nabla)$ is conformally-projectively flat.
\end{Coro}


\begin{thebibliography}{6}    

 
\bibitem{hase}I. Hasegawa , K. Yamauchi, Conformal-projective flatness of tangent bundle with complete lift statistical structure, Differential Geometry-Dynamical Systems, 10, 148-158, 2008.
\bibitem{kuro1}T. Kurose, On the divergence of 1-conformally flat statistical manifolds, T\^ohoku Math. J., 46, 427-433, 1994.  
\bibitem{kuro2}T. Kurose, Conformal-Projective geometry of Statistical Manifolds, Interdisciplinary information Sciences, 8, 1, 89-100, 2002.
\bibitem{mats}H. Matsuzoe, On realization of conformally-projectively flat statistical manifolds and the divergences, Hokkaido Math. J., 27, 409-421, 1998.
\bibitem{uoha}K. Uohashi, On $\alpha-$conformal equivalence of statistical manifolds, J. Geom., 75, 179-184, 2002.
\bibitem{zhang}J. Zhang, A note on curvature of $\alpha-$connections of a statistical manifold, AISM, 59, 161-170, 2007.


\end{thebibliography}
\end{document}